\documentclass[preprint,11pt]{elsarticle}

\usepackage{amssymb}
\usepackage{amsfonts}
\usepackage{}
  \usepackage{color}
 \usepackage{bbm}
 \usepackage{mathrsfs,latexsym}
 \usepackage[all]{xy}
 \usepackage{amsthm,amsmath,amsfonts,amssymb}
\newtheorem{thm}{Theorem}[section]
\newtheorem{lem}{Lemma}[section]

\newtheorem{cor}{Corollary}[section]
\newtheorem{pro}{Proposition}[section]
\newtheorem{rem}{Remark}[section]

\theoremstyle{definition}

\journal{}

\begin{document}
\numberwithin{equation}{section}
\renewcommand{\theequation}{\arabic{section}.\arabic{equation}}

\begin{frontmatter}

\title{Viscosity solutions to quaternionic Monge-Amp\`{e}re equations}

\author[label1]{Dongrui Wan}
\author[label2]{Wei Wang}
\address[label1]{College of Mathematics and Computational Science, Shenzhen University, Shenzhen, 518060, P. R. China, Email: wandongrui@szu.edu.cn}
\address[label2]{Department of Mathematics,
Zhejiang University, Zhejiang 310027,
 P. R. China, Email:   wwang@zju.edu.cn}

\begin{abstract}Quaternionic Monge-Amp\`{e}re equations have recently been studied intensively using methods from pluripotential theory. We present an alternative approach by using the viscosity methods. We study the viscosity solutions to the Dirichlet problem for quaternionic Monge-Amp\`{e}re equations $det(f)=F(q,f)$ with boundary value $f=g$ on $\partial\Omega$. Here $\Omega$ is a bounded domain on the quaternionic space $\mathbb{H}^n$, $g\in C(\partial\Omega)$, and $F(q,t)$ is a continuous function on $\Omega\times\mathbb{R}\rightarrow\mathbb{R}^+$ which is non-decreasing in the second variable. We prove a viscosity comparison principle and a solvability theorem. Moreover, the equivalence between viscosity and pluripotential solutions is showed.
\end{abstract}


\end{frontmatter}

\section{Introduction}
Viscosity methods provide an efficient tool for the study of weak solutions to nonlinear elliptic partial differential equations (see e.g \cite{caff,guide,ishii}, and references therein). The viscosity theory for real Monge-Amp\`{e}re equations has been developed by Ishii and Lions \cite{ishii}, and more generally for real $k$-Hessian equations by Trudinger and Wang \cite{hessian2,hessian3}. The viscosity approach for the complex case hasn't been studied until recent years (see e.g. \cite{riemann,harvey09,harvey11}). Viscosity solutions to degenerate complex Monge-Amp\`{e}re equations were studied on compact k\"{a}hler manifolds \cite{guedj}. Then Wang \cite{Y.Wang} presented a different viscosity approach to the Dirichlet problem for the complex Monge-Amp\`{e}re equations. More generally, viscosity solutions to complex Hessian equations have been studied in \cite{lu}. There is a nice survey on these last developments (see \cite{Zeriahi}).

From recent developments on viscosity theory for complex Monge-Amp\`{e}re equations, it is nature to develop such a treatment for quaternionic Monge-Amp\`{e}re equations, which is the main purpose of this paper. To the best of our knowledge, there is no reference on viscosity solutions to quaternionic Monge-Amp\`{e}re equations.

Quaternionic version of Calabi-Yau conjecture, which has applications in superstring theory, has attracted many analysts to study on this problem. It has already been solved on compact manifolds with a flat hyperK\"{a}hler metric \cite{alesker9}. It is meaningful and interesting to develop pluripotential theory on quaternionic manifold \cite{alesker2,alesker6,wang1}.

The quaternionic Monge-Amp\`{e}re operator is defined as the Moore determinant of the quaternionic Hessian of $u$:
 \begin{equation*}det(u)=det\left[\frac{\partial^2u}{\partial q_j\partial \bar{q}_k}(q)\right].\end{equation*}Alesker extended the definition of quaternionic Monge-Amp\`{e}re operator to continuous quaternionic plurisubharmonic functions and developed a basic theory of quaternionic plurisubharmonic functions in \cite{alesker1}. To define the quaternionic Monge-Amp\`{e}re operator on general quaternionic manifolds, Alesker introduced in \cite{alesker2} an operator in terms of the Baston operator $\triangle$, which is the first operator of the quaternionic complex on quaternionic manifolds. The $n$-th power of this operator is exactly the quaternionic Monge-Amp\`{e}re operator when the manifold is flat.

 In \cite{wan3} we established an explicit expression of the quaternionic Monge-Amp\`{e}re operator by using the operator $d_0$ and $d_1$, which are the first-order differential operators acting on the quaternionic versions of differential forms on the flat quaternionic space $\mathbb{H}^n$. We wrote down explicitly the quaternionic Monge-Amp\`{e}re operator in terms of real variables. This definition is much more convenient than using the Moore determinant. By defining the notion of closed positive currents in the quaternionic case, we generalized Bedford-Taylor theory, i.e., extend the definition of the quaternionic Monge-Amp\`{e}re operator to locally bounded quaternionic plurisubharmonic functions and prove the corresponding convergence theorem. As an application, the first author and Zhang proved in \cite{wan4} the quasicontinuity theorem and a comparison principle for quaternionic plurisubharmonic functions in $\mathbb{H}^n$.

The results about the quaternionic Monge-Amp\`{e}re operator proved by the first author \cite{wan4} and the second author \cite{wan3,wang1} play key role in this paper. This is the reason why we can get a viscosity approach to the quaternionic Monge-Amp\`{e}re equations on the quaternionic space $\mathbb{H}^n$.

Specifically, in this paper we consider the viscosity solutions to the following quaternionic Monge-Amp\`{e}re equation: \begin{equation}\label{1} -det(f)+F(q,f)=0,\end{equation}
with boundary value $f=g$ on $\partial\Omega$. Here $\Omega$ is a bounded domain on the quaternionic space $\mathbb{H}^n$, $g\in C(\partial\Omega)$, and $F(q,t):\Omega\times\mathbb{R}\rightarrow\mathbb{R}^+$ is a continuous function which is non-decreasing in the second variable.

 Our main results are the following:
\begin{thm}\label{t3} Assume that there exist a bounded subsolution $u$ and a bounded supersolution $v$ to (\ref{1}) such that $u_*=v^*=g$ on $\partial\Omega$. Then there exists a unique viscosity solution to (\ref{1}) with boundary value $g$. It is also the unique pluripotential solution. (Here denote by $u_*$ the lower semicontinuous envelope of $u$ and $v^*$ the upper semicontinuous envelope of $v$.)
\end{thm}
Moreover, we clarify the connection between viscosity solutions in our sense and pluripotential solutions.
\begin{thm}\label{t1}  Let $f$ be a bounded upper semi-continuous function in $\Omega$. Then the inequality \begin{equation}\label{2}det (f)\geq F(q,f)\end{equation} holds in the viscosity sense if and only if $f$ is plurisubharmonic in $\Omega$ and (\ref{2}) holds in the pluripotential sense.
\end{thm}

A key ingredient for our approach is the following viscosity comparison principle, which implies uniqueness of viscosity solutions to the quaternionic Monge-Amp\`{e}re equation (\ref{1}).
\begin{thm}\label{t2} (viscosity comparison principle) Let $u$ be a bounded viscosity subsolution and $v$ be a bounded viscosity supersolution of (\ref{1}). If $u\leq v$ on $\partial\Omega$ then $u\leq v$ on $\Omega$.
\end{thm}

\section{Preliminaries}
In this section, first we introduce the following definitions and basic results of quaternionic Monge-Amp\`{e}re operator and plurisubharmonic function.

For a point $q=(q_0,\ldots,q_{n-1})^t\in \mathbb{H}^n$, where $ ^t$ is the transpose, write
\begin{equation}\label{q}q_j=x_{4j}+x_{4j+1}\textbf{i}+x_{4j+2}\textbf{j}+x_{4j+3}\textbf{k},\end{equation} $j=0,\ldots,n-1.$ The Cauchy-Fueter operator is defined as
$$\overline{\partial}f=\left(
                         \begin{array}{c}
                           \overline{\partial}_{q_0}f \\
                           \vdots \\
                           \overline{\partial}_{q_{n-1}}f \\
                         \end{array}
                       \right)$$ for a $\mathbb{H}$-valued $C^1$ function $f$, where $$\overline{\partial}_{q_j}=\partial_{x_{4j}}+\textbf{i}\partial_{x_{4j+1}}+\textbf{j}\partial_{x_{4j+2}}+\textbf{k}\partial_{x_{4j+3}}.$$ Set $$\partial_{q_j}f:=\overline{\overline{\partial}_{q_j}\bar{f}}=\partial_{x_{4j}}f-\partial_{x_{4j+1}}f\textbf{i}-\partial_{x_{4j+2}}f\textbf{j}-\partial_{x_{4j+3}}f\textbf{k}.$$

  A real valued function $f:\mathbb{H}^n\rightarrow \mathbb{R}$ is called \emph{quaternionic plurisubharmonic} (PSH, for short) if it is upper semi-continuous and its restriction to any right quaternionic line is subharmonic (in the usual sense). Any quaternionic plurisubharmonic function is subharmonic.

  A square quaternionic matrix $A=(a_{jk})$ is called \emph{hyperhermitian} if its quaternionic conjugate $A^*=A$, or explicitly $a_{jk}=\overline{a_{kj}}$. The Moore determinant denoted by $det$ is defined on the class of all hyperhermitian matrices and takes real value.

\begin{pro}\label{p1.1}
$(1)$ (Proposition 2.1.6 in \cite{alesker1}) A real valued $C^2$ function $f$ on the domain $\Omega\subseteq\mathbb{H}^n$ is quaternionic plurisubharmonic if and only if at every point $q\in\Omega$ the hyperhermitian matrix $$\left[\frac{\partial^2f}{\partial q_j\partial \bar{q}_k}(q)\right],$$ called the \emph{quaternionic Hessian}, is non-negative definite.\\
$(2)$ (Proposition 2.1 in \cite{wang1}) For a hyperhermitian $(n\times n)$-matrix $A$, there exists a unitary matrix $U$ such that $U^*AU$ is diagonal and real.\\
$(3)$ (Corollary 2.1 in \cite{wang1}) For a  real valued $C^2$ function $u$, the matrix $\left[\frac{\partial^2u}{\partial q_j\partial \bar{q}_k}\right]$ is hyperhermitian. We have $\left[\frac{\partial^2u}{\partial \bar{q}_j\partial q_k}(q)\right]=\left[\frac{\partial^2u}{\partial q_j\partial \bar{q}_k}(q)\right]^t,$ which is also a positive hyperhermitian matrix if $u$ is plurisubharmonic. And $det\left(\frac{\partial^2u}{\partial \bar{q}_j\partial q_k}(q)\right)=det\left(\frac{\partial^2u}{\partial q_j\partial \bar{q}_k}(q)\right)$ . We will denote for brevity $det(u):=det\left(\frac{\partial^2u}{\partial q_j\partial \bar{q}_k}(q)\right)$.\\
$(4)$ (Lemma 2.1 in \cite{wang1}) For any quaternionic $(n\times n)$-matrix $X$, let $X^{\mathbb{R}}$ be the \emph{real $(4n\times 4n)$-matrix associated to $X$}, which was introduced by the second author \cite{wang1}. $X^{\mathbb{R}}$ satisfies: $$(Xq)^{\mathbb{R}}=X^{\mathbb{R}}q^{\mathbb{R}},$$ where $q^{\mathbb{R}}$ is the associated real vector in $\mathbb{R}^{4n}$ defined as $$q^{\mathbb{R}}:=(x_0,x_1,\ldots,  x_{4n-1})^t,$$ for $q\in \mathbb{H}^n$ given by (\ref{q}).
$X^{\mathbb{R}}$ is symmetric if $X$ is hyperhermitian. $X^{\mathbb{R}}$ is positive if $X$ is positive hyperhermitian. \\
$(5)$ (Theorem 1.1.17 in \cite{alesker1}) The function $X\mapsto[det(X)]^{\frac{1}{n}}$ is concave on the cone of the positive definite hyperhermitian matrices. If $A,B\geq0$, then $det(A+B)\geq det(A)+det(B)$.
\end{pro}

In \cite{wan3} we established an explicit expression of the quaternionic Monge-Amp\`{e}re operator by using the Baston operator $\triangle$. We wrote down explicitly the quaternionic Monge-Amp\`{e}re operator in terms of real variables. By defining the notion of quaternionic closed positive currents, we proved the following results which are necessary in this paper.

Let $u_1,\ldots,u_n$ be locally bounded PSH functions in $\Omega$. Then $\triangle u_1\wedge\ldots\wedge\triangle u_n$ can be defined inductively as a positive Radon measure (see Proposition 3.9 in \cite{wan3}).
\begin{pro}\label{p2.3} (1) (Theorem 3.1 in \cite{wan3}) Let $v^1,\ldots,v^n\in PSH\cap L_{loc}^{\infty}(\Omega)$. Let $\{v_j^1\}_{j\in\mathbb{N}},\ldots,\{v_j^n\}_{j\in\mathbb{N}}$ be decreasing sequences of $PSH$ functions in $\Omega$ such that $\lim_{j\rightarrow\infty}v_j^t=v^t$   pointwisely in $\Omega$ for each $t $. Then the currents $\triangle v_j^1\wedge\ldots\wedge\triangle v_j^n $ converge weakly to $\triangle v^1\wedge\ldots\wedge\triangle v^n $ as $j\rightarrow\infty$.\\
(2) (Proposition 4.1 in \cite{wan4}) Let $\{u_j\}_{j\in \mathbb{N}}$ be a sequence in $PSH\cap L_{loc}^\infty(\Omega)$ that increases to $u\in PSH\cap L_{loc}^\infty(\Omega)$ almost everywhere in $\Omega$ $($with respect to Lebesgue measure$)$. Then the currents $(\triangle u_j)^n$ converge weakly to $(\triangle u)^n$ as $j\rightarrow\infty$.
\end{pro}

\begin{pro}\label{l2.10}(pluripotential comparison principle, Corollary 1.1 in \cite{wan4}) Let $u,v\in PSH\cap L_{loc}^\infty(\Omega)$. If for any $\zeta\in \partial\Omega$, $$\liminf_{\zeta\leftarrow q\in \Omega}(u(q)-v(q))\geq0,$$ and $(\triangle u)^n\leq (\triangle v)^n$ in $\Omega$, then $u\geq v$ in $\Omega$.
\end{pro}
\begin{lem}\label{t2.1}(Appendix A. in \cite{wan3}) Let $f_1,\ldots,f_n$ be $C^2$ functions in
$\mathbb{H}^n$. Then we have
\begin{equation}\label{2.17}\triangle_n(f_1,\ldots,
f_n)=n!~\text{det}~(f_1,\ldots,f_n),
\end{equation}where $\text{det}~(f_1,\ldots,f_n)$ is the mixed Monge-Amp\`{e}re operator and $\triangle_n(f_1,\ldots, f_n)$ is the
coefficient of the $2n$-form $\triangle f_1\wedge\ldots\wedge\triangle
f_n$.
\end{lem}

\begin{rem}As in \cite{wan3} it follows from Lemma \ref{t2.1} that the operator $f\mapsto\triangle_nf$
coincides with the quaternionic Monge-Amp\`{e}re operator $f\mapsto det(f)$. Proposition \ref{p2.3} and Proposition \ref{l2.10} also hold for the quaternionic Monge-Amp\`{e}re operator $det(f)$.
\end{rem}

See \cite{alesker1,wan3,wan4,wang1} for more information about quaternionic  Monge-Amp\`{e}re operator and quaternionic closed positive currents.

\section{viscosity solution and pluripotential solution}
In this section, we introduce the definitions of viscosity subsolutions and viscosity supersolutions to the quaternionic Monge-Amp\`{e}re equations, and clarify their relation with pluripotential ones.

Let $\Omega$ be a bounded domain in $\mathbb{H}^n$. Let $f:\Omega\rightarrow \mathbb{R}\cup\{-\infty\}$ be a function and let $\varphi$ be a $C^2$ function in a neighborhood of $q_0\in \Omega$. We say that \emph{$\varphi$ touches $f$ from above (resp. from below) at $q_0$} if $\varphi(q_0)=f(q_0)$ and $\varphi(q)\geq f(q)$ (resp. $\varphi(q)\leq f(q)$) for every $q$ in a neighborhood of $q_0$.

An upper semicontinuous function $f:\Omega\rightarrow\mathbb{R}\cup\{-\infty\}$ is a \emph{viscosity subsolution} to (\ref{1})
if $f\not\equiv -\infty$ and for any $q_0\in \Omega$ and any $C^2$ function $\varphi$ which touches $f$ from above at $q_0$ then $$det (\varphi)\geq F(q,\varphi)~~~\text{at}~q_0.$$
We also say that $ det (f)\geq F(q,f)$ holds\emph{ in the viscosity sense}.

A lower semicontinuous function $f:\Omega\rightarrow\mathbb{R}\cup\{+\infty\}$ is a \emph{viscosity supersolution} to (\ref{1})
if $f\not\equiv +\infty$ and for any $q_0\in \Omega$ and any $C^2$ function $\varphi$ which touches $f$ from below at $q_0$ then $$det (\varphi)_+\leq F(q,\varphi)~~~\text{at}~q_0.$$
Here $det (\varphi)_+$ is defined to be itself if $\varphi$ is plurisubharmonic and $0$ otherwise.

If $u\in C^2(\Omega)$ then $det (u)\geq F(q,u)$ (or $det (u)_+\leq F(q,u)$) holds in the viscosity sense iff it holds in the usual sense.

A function $f:\Omega\rightarrow\mathbb{R}$ is a \emph{viscosity solution} to (\ref{1}) if it is both a viscosity subsolution and a viscosity  supersolution to (\ref{1}). A viscosity solution is automatically a continuous function in $\Omega$.

Note that the class of viscosity subsolutions is stable under taking maximum. It is also stable along monotone sequences as the following lemma shows.
\begin{lem}\label{l2.1} Let $\{f_j\}$ be a monotone sequence of viscosity subsolutions of (\ref{1}). If $f_j$ is uniformly bounded from above and $f:=(\lim_j f_j)^*\not\equiv-\infty$ then $f$ is also a viscosity subsolution of (\ref{1}).
\end{lem}
\proof This proof can be found in \cite{guide}. For convenience, we repeat it here. Take $q_0\in \Omega$ and a $C^2$ function $\varphi$ in $B(q_0,r)\subseteq\Omega$ which touches $f$ from above at $q_0$. We can choose a sequence $\{q_j\}$ in $ B=B(q_0,\frac{r}{2})$ converging to $q_0$ and a subsequence of $\{f_j\}$ (denoted also by $\{f_j\}$) such that $f_j(q_j)\rightarrow f(q_0)$. Fix $\varepsilon>0$. For each $j$, let $y_j$ be the maximum point of $f_j-\varphi-\varepsilon|q-q_0|^2$ in $B$. So \begin{equation}\label{2.1}f_j(q_j)-\varphi(q_j)-\varepsilon|q_j-q_0|^2\leq f_j(y_j)-\varphi(y_j)-\varepsilon|y_j-q_0|^2.\end{equation}
Assume that $y_j\rightarrow y\in B$. Letting $j\rightarrow+\infty$ in (\ref{2.1}) and noting that $\limsup f_j(y_j)\leq f(y)$, we have $$0\leq f(y)-\varphi(y)-\varepsilon|y-q_0|^2.$$ Since $\varphi$ touches $f$ from above at $q_0$, we get that $y=q_0$. So $y_j\rightarrow q_0$ as $j\rightarrow+\infty$. Then by (\ref{2.1}) again we have $f_j(y_j)\rightarrow f(q_0)$. For $j$ large sufficiently, the function $\varphi(q)+\varepsilon|q-q_0|^2+ f_j(y_j)-\varphi(y_j)-\varepsilon|y_j-q_0|^2$ touches $f_j$ from above at $y_j$. Since $f_j$ is a viscosity subsolution of (\ref{1}), by definition we have $$det(\varphi+\varepsilon|q-q_0|^2)(y_j)\geq F(y_j,f_j(y_j)).$$
Let $j\rightarrow+\infty$ to get $$det(\varphi+\varepsilon|q-q_0|^2)(q_0)\geq F(q_0,f(q_0)).$$
Hence $f$ is also a viscosity subsolution of (\ref{1}).
\endproof

When $F\equiv0$, the viscosity subsolutions of (\ref{1}) are precisely the PSH functions on $\Omega$.
\begin{pro}\label{p2.1}A function $f$ is plurisubharmonic on $\Omega$  if and only if it is a viscosity subsolution of $-det( f)=0$.
\end{pro}
To prove this proposition, we need the following lemmas first.

\begin{lem}\label{l2.2} Let $Q$ be a hyperhermitian matrix such that for any semipositive hyperhermitian matrix $H$, $det(Q+H)\geq0$. Then $Q$ is semipositive.
\end{lem}This lemma can be proved by diagonalizing $Q$ by Proposition \ref{p1.1} (2).

The second author introduced in \cite{wang1} the operator  \begin{equation}\label{2.2}\Delta_\mathbbm{a} v:=\frac{1}{2}Re \sum_{j,k=1}^n a_{kj}\frac{\partial^2v}{\partial \bar{q}_j\partial q_k}=\frac{1}{2}ReTr\left(\mathbbm{a}(\frac{\partial^2v}{\partial \bar{q}_j\partial q_k})\right),
\end{equation}for a quaternionic positive hyperhermitian $(n\times n)$-matrix $\mathbbm{a}=(a_{jk})$ and a $C^2$ real function $v$. This is an elliptic operator of constant coefficients. This operator is the quaternionic counterpart of complex K\"{a}hler operator, which plays key role in the viscosity approach for the complex case.
\begin{lem}\label{l2.3}(Lemma 4.1 in \cite{wang1}) If $v$ is subharmonic with respect to $\Delta_\mathbbm{a}$ for any quaternionic positive hyperhermitian $(n\times n)$-matrix $\mathbbm{a}$, then $v$ is plurisubharmonic.
\end{lem}
With the help of this operator, now we can prove Proposition \ref{p2.1}. We take idea from \cite{guedj,lu}.\vspace{2mm}

\emph{Proof of Proposition \ref{p2.1}} Assume that $f\in PSH(\Omega)$. Denote by $f_\varepsilon$ its standard smooth regularization, then $f_\varepsilon$ is also PSH and smooth. Hence $det( f_\varepsilon)\geq0$, i.e., $f_\varepsilon$ is a classical subsolution of $-det( f)=0$. Since $f_\varepsilon$ converges decreasingly to $f$ as $\varepsilon\rightarrow0$, it follows from Lemma \ref{l2.1} that $f$ is a viscosity subsolution of $-det( f)=0$.

Conversely, assume that $f$ is a viscosity subsolution of $-det( f)=0$. Let $q_0$ be such that $f(q_0)\neq-\infty$ and let $\varphi\in C^2(\{q_0\})$ be such that $\varphi$ touches $f$ from above at $q_0$. Then for any $h(q-q_0)$ whose quaternionic Hessian matrix $H(q-q_0)$ is hyperhermitian semipositive, $\varphi+h(q-q_0)$ also touches $f$ from above at $q_0$. By the definition of viscosity subsolution we have $det(\varphi+h(q-q_0)) \geq0$. Since the quaternionic matrix $\left[\frac{\partial^2\varphi}{\partial \bar{q}_j\partial q_k}\right]$ is hyperhermitian, it is semipositive by Lemma \ref{l2.2}.

It follows from (\ref{2.2}) and the semipositivity of $[\frac{\partial^2\varphi}{\partial \bar{q}_j\partial q_k}]$ that $\Delta_\mathbbm{a}\varphi\geq0$ for any quaternionic positive hyperhermitian $(n\times n)$-matrix $\mathbbm{a}$. Therefore $\Delta_\mathbbm{a} f\geq0$ in the viscosity sense. In appropriate coordinates this constant coefficient operator is just the Laplace operator. Then Prop.3.2.10' in \cite{homander} implies that $f$ is $\Delta_\mathbbm{a}$-subharmonic, and so is $L_{loc}^1$ and satisfies $\Delta_\mathbbm{a}f\geq0$ in the sense of distributions. It follows from Lemma \ref{l2.3} that $f$ is PSH.
 \qed

\begin{cor}Lemma \ref{l2.1} still holds if the sequence $\{f_j\}$ is not monotone.
\end{cor}
\proof For each $j$, let $$u_j=(\sup_{k\geq j}f_k)^*,v_l=\max\{f_j,\ldots,f_{j+l}\}.$$ Since the class of viscosity subsolutions is stable under taking maximum, $v_l$ is a viscosity subsolution of (\ref{1}). Note that $u_j=(\sup_{l\geq0}v_l)^*$ and the sequence $\{v_l\}$ is monotone. It follows from Lemma \ref{l2.1} that $u_j$ is a viscosity subsolution of (\ref{1}). By Proposition \ref{p2.1}, each $f_j$ and $u_j$ are plurisubharmonic. Note that $u_j$ converges decreasingly to $f$ as $j\rightarrow+\infty$. The proof is complete.
\endproof

Recall that when $f$ is plurisubharmonic and locally bounded, its quaternionic Monge-Amp\`{e}re measure $det( f)$ is well defined and is continuous on decreasing sequences (as showed in Section 2 and in \cite{wan3}). Now let us make the basic connection between this pluripotential notion and its viscosity counterpart when the function $F(q,t)$ in (\ref{1}) does not depend on $t$.

\begin{lem} \label{l2.4}(Proposition 2.6 in \cite{wang1}) For a positive quaternionic hyperhermitian $(n\times n)$-matrix $X$, $$(detX)^{\frac{1}{n}}=\frac{1}{n}\inf_\mathbbm{a} ReTr(\mathbbm{a} X)=\frac{1}{n}\inf_\mathbbm{a} Tr(\mathbbm{a} X\mathbbm{a}),$$ where the infimums are taken over all positive quaternionic hyperhermitian $(n\times n)$-matrices $\mathbbm{a}$ such that $det \mathbbm{a}\geq1$.
\end{lem}
The lemma above is the key ingredient in the proof of the following result. We take idea from \cite{guedj,lu,Y.Wang}.
\begin{pro}\label{p2.2}Let $f$ be a bounded upper semicontinuous function in $\Omega$ and let $g\geq0$ be a continuous function.\\
(1). If $f$ is plurisubharmonic such that
\begin{equation}\label{2.3}det( f)\geq g
\end{equation}in the pluripotential sense then it also holds in the viscosity sense.\\
(2). Conversely, if (\ref{2.3}) holds in the viscosity sense then $f$ is plurisubharmonic and the inequality holds in the pluripotential sense.
\end{pro}
\proof Let $f\in PSH\cap L_{loc}^\infty(\Omega)$ satisfy (\ref{2.3}). Assume that a $C^2$ function $\varphi$ touches $f$ from above at $q_0\in \Omega$. Suppose by contradiction that $det( \varphi)(q_0)<g(q_0)$. By choosing $\varepsilon>0$ small enough and letting $\varphi_\varepsilon:=\varphi+\varepsilon|q-q_0|^2$, we have $0<det (\varphi_\varepsilon)<g$ in a neighborhood of $q_0$ by the continuity of $g$. It follows from the proof of Proposition \ref{p2.1} that $\varphi_\varepsilon$ is plurisubharmonic in a neighborhood of $q_0$, say $B$. Now for $\delta>0$ small enough, we have $\varphi_\varepsilon-\delta\geq f$ near $\partial B$ and $det(\varphi_\varepsilon)\leq det( f)$. The pluripotential comparison principle (Proposition \ref{l2.10}) yields $\varphi_\varepsilon-\delta\geq f$ on $B$. But this fails at $q_0$. Hence $det(\varphi)(q_0)\geq g(q_0)$ and (\ref{2.3}) holds in the viscosity sense.

To prove (2), we first assume that $g>0$ is smooth. Suppose that $\varphi\in C^2$ touches $f$ from above at $q_0\in \Omega$. Then $det( \varphi)\geq g$ at $q_0$, which means the hyperhermitian matrix $Q=\left[\frac{\partial^2\varphi}{\partial \bar{q}_j\partial q_k}(q_0)\right]$ satisfies $det(Q)\geq g>0$. By Lemma \ref{l2.4} and (\ref{2.2}), we have
\begin{equation}\label{*}(det Q)^\frac{1}{n}=\frac{2}{n}\inf_\mathbbm{a} \Delta_\mathbbm{a}\varphi\geq g^{\frac{1}{n}},\end{equation}
for every positive hyperhermitian $(n\times n)$-matrix $\mathbbm{a}=(a_{jk})$ with $det \mathbbm{a}\geq1$. Hence $\Delta_\mathbbm{a}\varphi\geq \frac{n}{2}g^{\frac{1}{n}},$ i.e., $f$ is a viscosity subsolution of the constant coefficient elliptic equation $\Delta_\mathbbm{a}\varphi= \frac{n}{2}g^{\frac{1}{n}}.$ Choose a $C^2$ solution $h$ of $\Delta_\mathbbm{a}\varphi= \frac{n}{2}g^{\frac{1}{n}}$ in a neighborhood of $q_0$. Then $u:=f-h$ satisfies $\Delta_\mathbbm{a} u\geq0$ in the viscosity sense. Once again, Proposition 3.2.10' in \cite{homander} implies that $u$ is $\Delta_\mathbbm{a}$-subharmonic and satisfies $\Delta_\mathbbm{a} u\geq0$ in the sense of distributions. Hence $\Delta_\mathbbm{a} f\geq \frac{n}{2}g^{\frac{1}{n}}$ in the sense of distributions. Consider the standard smooth regularization $f_\varepsilon:=f*\rho_\varepsilon$. We see that $ \Delta_\mathbbm{a} f_\varepsilon\geq \frac{n}{2}(g^{\frac{1}{n}})_\varepsilon$. By (\ref{*}) we obtain $(det(f_\varepsilon))^{\frac{1}{n}}\geq (g^{\frac{1}{n}})_\varepsilon$. Thus $det (f_\varepsilon)\geq ((g^{\frac{1}{n}})_\varepsilon)^n .$ Letting $\varepsilon\rightarrow0$ and noting that the quaternionic Monge-Amp\`{e}re operator is continuous on decreasing sequence (Proposition \ref{p2.3}), we get $det( f)\geq g$.

In case $g>0$ is only continuous we observe that $$g=\sup\{h\in C^\infty(\Omega),0<h\leq g\}.$$ If $det( f)\geq g$ in the viscosity sense then we also have $det( f)\geq h$. By the proof above, $det( f)\geq h$ in the pluripotential sense for every $0<h\leq g$. We conclude that $det( f)\geq g$ holds in the pluripotential sense.

Now let $g\geq 0$ be continuous. Consider $f_\delta(q):=f(q)+\delta|q|^2$. Since $det(|q|^2)=8^n$, $det( f_\delta)\geq (g+8^n\delta^n)$ in the viscosity sense. By the proof above we have $det( f_\delta)\geq (g+8^n\delta^n)$ in the pluripotential sense. The result follows by letting $\delta\rightarrow0$.
\endproof

\emph{Proof of Theorem \ref{t1}} Assume that (\ref{2}) holds in the viscosity sense. We approximate $f$ by its sup-convolution:
\begin{equation}\label{3.4}f^\delta(q)=\sup_{q'\in \Omega}\{f(q')-\frac{1}{2\delta^2}|q-q'|^2\},~~~q\in \Omega_\delta,\end{equation}for $\delta>0$ small enough, where $\Omega_\delta:=\{q\in\Omega:dist(q,\partial\Omega)>A\delta\}$ and $A>0$ is a large constant such that $A^2>2 osc_\Omega f$. This family of continuous semi-convex functions decreases towards $f$ as $\delta\rightarrow0$. As in \cite{guedj,unique}, $f^\delta$ satisfies the following inequality in the viscosity sense:
\begin{equation}\label{2.5}det( f^\delta)\geq F_\delta(q,f^\delta)~~~~~\text{in}~~\Omega_\delta, \end{equation}
where $F_\delta(q,t)=\inf_{|q'-q|\leq A\delta} F(q',t)$. It follows from Proposition \ref{p2.1} that $f^\delta$ is plurisubharmonic. Apply Proposition \ref{p2.2} to get that (\ref{2.5}) holds in the pluripotential sense. Since $F$ is non-decreasing in the second variable, we have $$det( f^\delta)\geq F_\delta(q,f^\delta)\geq F_\delta(q,f).$$ Note that $F_\delta(q,t)$ increases towards $F(q,t)$ as $\delta\rightarrow0$. Since the quaternionic  Monge-Amp\`{e}re operator is continuous on decreasing sequence of plurisubharmonic functions (Proposition \ref{p2.3}), we finally obtain $det( f)\geq F(q,f)$ in the pluripotential sense.

Now we prove the other implication. Let $f$ be a PSH function such that (\ref{2}) holds in the pluripotential sense. If $f$ were continuous, then we could use Proposition \ref{p2.2}. For non-continuous $f$, we approximate $f$ by its sup-convolution $f^\delta$ as above. As in \cite{guedj,unique}, by Lemma \ref{l2.5} below we can get that (\ref{2.5}) holds in $\Omega_\delta$ in the pluripotential sense. Apply Proposition \ref{p2.2} to get that $f^\delta$ is a viscosity subsolution of the equation $det( f)=F_\delta(q,f)$ in $\Omega_\delta$. We want to prove that $f$ is a viscosity subsolution of the equation $det( f)=F(q,f)$. Now it suffices to let $\delta\rightarrow0$.

\qed

\begin{lem}\label{l2.5} Let $f$ be a bounded plurisubharmonic function in $\Omega$ such that $det( f)\geq F(q,f)$ in the pluripotential sense in $\Omega$, where $F:\Omega\times \mathbb{R}\rightarrow\mathbb{R}^+$ is a continuous function and is non-decreasing in the second variable. Then the sup-convolution $\{f^\delta\}$ given by (\ref{3.4}) satisfies $$det( f^\delta)\geq F_\delta(q,f^\delta)$$ in the pluripotential sense, where $$F_\delta(q,t):=\inf_{q'}\{F(q',t):|q-q'|\leq A\delta\}.$$
\end{lem}
\proof Fix $\delta>0$ small enough. For $q'\in B(0, A\delta)$, let $$f_{q'}(q):=f(q-q')-\frac{1}{2\delta^2}|q'|^2, ~~~q\in \Omega_\delta. $$ Observe that $f_{q'}$ is also a bounded PSH function on $\Omega_\delta$. Since the quaternionic Monge-Amp\`{e}re operator is invariant under translation transformation (proved by Alesker \cite{alesker4}), we get $det( f_{q'})\geq F_\delta(q,f_{q'})$ in the pluripotential sense on $\Omega_\delta$.

 Note that $f$ is the upper envelope of the family $\{f_{q'}:|q'|\leq A\delta\}$. It follows from the well known Choquet lemma that there exists a sequence of points $\{q'_j\}_{j\in \mathbb{N}}$ with $|q_j'|\leq A\delta$ such that $f^\delta=(\sup_j f_{q'_j})^*$ on $\Omega_\delta$. Let $\theta_j:=\sup_{0\leq k\leq j}f_{q'_k}$. Then $\{\theta_j\}$ is an increasing sequence of bounded plurisubharmonic functions that converge almost everywhere to $f^\delta$ on $\Omega_\delta$. We claim that $\theta_j$ also satisfies the inequality $det( \theta_j)\geq F_\delta(q,\theta_j)$ in the pluripotential sense. By Demailly's inequality for quaternionic Monge-Amp\`{e}re operator we proved in Proposition 3.6 in \cite{wan5},
$$det( \max\{u_1,u_2\})\geq \chi_{\{u_1>u_2\}}det(u_1)+\chi_{\{u_1\leq u_2\}}det( u_2).$$
If $det( u_i)\geq F_\delta(q,u_i)$ for $i=1,2$, then the function $\omega:=\max\{u_1,u_2\}$ satisfies the same inequality, i.e. $det( \omega)\geq F_\delta(q,\omega)$. This proves our claim that $det( \theta_j)\geq F_\delta(q,\theta_j)$. Since $\{\theta_j\}$ converges increasingly almost everywhere in $\Omega_\delta$ to $f^\delta$,  we get $det( f^\delta)\geq F_\delta(q,f^\delta)$ in the pluripotential sense by Proposition \ref{p2.3} (2).
\endproof

\section{proof of the main results}
In this section, first we prove an important viscosity comparison principle which plays a key role in the proof of our main theorem. We use the classical method from \cite{caff,guide,guedj,lu,Y.Wang}.

The following results are a theorem of Alexandroff-Buselman-Feller (see Theorem 1 in \cite{measure}, or Section 1.2, Appendix 2 in \cite{nonlinear}) and the ABP estimate (see Theorem  5 and Corollary 2 in \cite{Y.Wang}).

Since the quaternionic Monge-Amp\`{e}re operator can be defined in terms of real variables, as is shown in Section 2, we claim that these results also hold for real functions over the quaternionic space $\mathbb{H}^n$ by identifying $\mathbb{H}^n$ with $\mathbb{R}^{4n}$.

A real function $u$ on $\Omega\subseteq\mathbb{H}^n$ is called \emph{semi-concave} (resp. \emph{semi-convex}) if there exists $K>0$ (resp. $K<0$) such that for every $q_0\in \Omega$, there exists a quadratic polynomial $P=K|q|^2+l$, which touches $u$ from above (resp. below) at $q_0$. Here $l$ is an affine function.

A real function $u$ on $\Omega\subseteq\mathbb{H}^n$ is called \emph{punctually second order differentiable} at $q_0\in \Omega$, if there exists a quadratic polynomial $\varphi$ such that $u(q)=\varphi(q)+o(|q-q_0|^2)$ as $q\rightarrow q_0$.

\begin{lem}\label{l2.6}Every continuous semi-convex (or semi-concave) function is punctually second order differentiable almost everywhere.
\end{lem}Sets in $\mathbb{H}^n$ are measured according to the standard $4n$-Lebesgue measure.
\begin{lem}\label{l2.7}(ABP estimate) Let $\omega\in C(\overline{\Omega})$ be a semi-concave function in a bounded domain $\Omega\subseteq\mathbb{H}^n$ such that $\Omega\subseteq B_d\subseteq B_{2d}$. Assume that $\omega\geq0$ on $\partial\Omega$ and $\min_\Omega \omega=-a$ for $a>0$. Let $E\subseteq \Omega$ be a set such that the Lebesgue measure of $\Omega\backslash E$ is zero. Define $$\Gamma_\omega(q):=\sup\{\varphi(q):\varphi ~\text{is~convex~in~}B_{2d}, \varphi\leq \min\{\omega,0\}~\text{in}~\Omega\}, ~~q\in B_{2d}.$$ Then for each $\delta\in (0,\frac{a}{2d})$, there exists a point $q_0\in E$ such that: \\
(1)  $\omega(q_0)=\Gamma_\omega(q_0)<0$; \\
(2)  $\Gamma_\omega$ is punctually second order differentiable at $q_0$;\\
(3) $$(det_\mathbb{R}D^2\Gamma_\omega(q_0))^{\frac{1}{4n}}\geq \frac{\delta}{d}.$$ Here $D^2 \Gamma_\omega$ is the real Hessian matrix of $\Gamma_\omega$ and $det_\mathbb{R}$ is the real determinant.
\end{lem}
\begin{lem}\label{l4.1}(Proposition 2.5 in \cite{wang1})  For a quaternionic hyperhermitian $(n\times n)$-matrix $X$ and a $C^2$ real function $u$, we have $$ReTr\left(X (\frac{\partial^2u}{\partial\bar{q}_j\partial q_k})\right)=Tr \left(X^{\mathbb{R}}(\frac{\partial^2u}{\partial x_s\partial x_t})\right).$$
\end{lem}
Recall that for any positive real symmetric $(n\times n)$-matrix $X$, $$(det_\mathbb{R}X)^{\frac{1}{n}}=\frac{1}{n}\inf_{\mathbbm{b}\in\Lambda' } Tr(\mathbbm{b} X),$$ where the infimum $\mathbbm{b}\in\Lambda' $ is taken over all positive real symmetric $(n\times n)$-matrices $\mathbbm{b}$ such that $det_\mathbb{R}\mathbbm{b}\geq1$ (See for example \cite{Gaveau} for details). Denote by
$det_\mathbb{R}$ the usual real determinant of real matrix and denote by $det$ the Moore determinant of quaternionic matrix.
\begin{pro}\label{p3.2} For a function $u\in PSH\cap C^2$  , we have the inequality: $$\left(det(\frac{\partial^2u}{\partial\bar{q}_j\partial q_k})\right)^{\frac{1}{n}}\geq 4\left(det_\mathbb{R}(\frac{\partial^2u}{\partial x_s\partial x_t})\right)^{\frac{1}{4n}}.$$
\end{pro}
\proof For $u\in PSH\cap C^2$, the quaternionic matrix $\left[\frac{\partial^2u}{\partial\bar{q}_j\partial q_k}\right]$ is positive hyperhermitian by Proposition \ref{p1.1} (3). By Lemma \ref{l2.4} and Lemma \ref{l4.1}, we have
\begin{equation*}\begin{aligned}\left(det(\frac{\partial^2u}{\partial\bar{q}_j\partial q_k})\right)^{\frac{1}{n}}&=\frac{1}{n}\inf_{\mathbbm{a}\in \Lambda} ReTr\left(\mathbbm{a} (\frac{\partial^2u}{\partial\bar{q}_j\partial q_k})\right)=\frac{1}{n}\inf_{\mathbbm{a}\in \Lambda} Tr\left(\mathbbm{a}^\mathbb{R}(\frac{\partial^2u}{\partial x_s\partial x_t})\right )
\\&\geq 4\left[\frac{1}{4n}\inf_{\mathbbm{b}\in \Lambda'} Tr\left(\mathbbm{b}(\frac{\partial^2u}{\partial x_s\partial x_t})\right )\right]= 4\left(det_\mathbb{R}(\frac{\partial^2u}{\partial x_s\partial x_t})\right)^{\frac{1}{4n}}.
\end{aligned}\end{equation*}The inequality is because, the infimum $\inf_{\mathbbm{a}\in \Lambda}$ is taken over all positive hyperhermitian $(n\times n)$-matrices $\mathbbm{a}$ and the $\inf_{\mathbbm{b}\in \Lambda'} $ is taken over all positive real symmetric $(4n\times 4n)$-matrices $\mathbbm{b}$. By Proposition \ref{p1.1} (4), $\mathbbm{a}^\mathbb{R}$ is positive symmetric for positive hyperhermitian $\mathbbm{a}$.
\endproof
\emph{Proof of Theorem \ref{t2}}~Note that $F$ is non-decreasing in the second variable. By replacing $u,v$ by $u-\varepsilon,v+\varepsilon$ we can assume that $u<v$ in a small neighborhood of $\partial\Omega$. Suppose by contradiction that there exists $q_0\in \Omega$ such that $u(q_0)-v(q_0)=a>0$.  As in the proof of Theorem \ref{t1}, we consider the sup-convolution $u^\varepsilon$ of $u$ defined by (\ref{3.4}) and the inf-convolution $v_\varepsilon$ of $v$ defined by
$$v_\varepsilon(q)=\inf_{|q'-q|\leq \varepsilon}\left\{v(q')+\frac{1}{2\varepsilon^2}|q-q'|^2\right\},~~~q\in \Omega_\varepsilon,$$ for $\varepsilon>0$ small enough.
Then $u^\varepsilon$ is semi-convex and $v_\varepsilon$ is semi-concave. By Dini lemma, $\omega_\varepsilon:=v_\varepsilon-u^\varepsilon\geq 0$ near $\partial\Omega$, for $\varepsilon>0$ small enough. Fix some open subset $U\Subset \Omega$ such that $\omega_\varepsilon\geq0$ on $\Omega\backslash U$. Fix $\varepsilon>0$ small enough. Denote by $E_\varepsilon$ the set of all points in $U$ where $\omega_\varepsilon,u^\varepsilon,v_\varepsilon$ are punctually second order differentiable. Then by Lemma \ref{l2.6} the Lebesgue measure of $U\backslash E_\varepsilon$ must be zero.

Fix some $d>0$ such that $\Omega\subseteq B_d\subseteq B_{2d}$. Define $$G_\varepsilon(q)=\sup\{\varphi(q):~\varphi \text{~is~convex~in~}B_{2d},~\varphi\leq\min(\omega_\varepsilon,0) ~\text{in}~\Omega\},~~q\in B_{2d}.$$
Since $\omega_\varepsilon\geq0$ on $\partial\Omega$ and $\omega_\varepsilon(q_0)=-a<0$, by Lemma \ref{l2.7} we can find $q_\varepsilon\in E_\varepsilon$ such that: $\omega_\varepsilon(q_\varepsilon)=G_\varepsilon(q_\varepsilon)<0$; $G_\varepsilon$ is punctually second order differential at $q_\varepsilon$ and $det_\mathbb{R}D^2G_\varepsilon(q_\varepsilon)\geq \delta$, where $\delta>0$ does not depend on $\varepsilon$. It follows from Proposition \ref{p3.2} that $det G_\varepsilon(q_\varepsilon)\geq \delta_1$, where $\delta_1>0$ does not depend on $\varepsilon$.

On the other hand, we have
$$det( u^\varepsilon)(q_\varepsilon)\geq F_\varepsilon(q_\varepsilon,u^\varepsilon(q_\varepsilon)),$$ as in the proof of Theorem \ref{t1}.
Since $(G_\varepsilon+u^\varepsilon)(q_\varepsilon)=v_\varepsilon(q_\varepsilon)$, $G_\varepsilon+u^\varepsilon$ touches $v_\varepsilon$ from below at $q_\varepsilon$. Note that $F$ is non-decreasing in the second variable and $\omega_\varepsilon(q_\varepsilon)<0$. Since $G_\varepsilon+u^\varepsilon$ is punctually second order differentiable at $q_\varepsilon$, we have
$$det (G_\varepsilon+u^\varepsilon)(q_\varepsilon)\leq F^\varepsilon(q_\varepsilon,v_\varepsilon(q_\varepsilon))\leq F^\varepsilon(q_\varepsilon,u^\varepsilon(q_\varepsilon)).$$
By Proposition \ref{p1.1} (5), we have $det (G_\varepsilon+u^\varepsilon)\geq det( G_\varepsilon)+det ( u^\varepsilon)$. Therefore
$$\delta_2+F_\varepsilon(q_\varepsilon,u^\varepsilon(q_\varepsilon))\leq F^\varepsilon(q_\varepsilon,u^\varepsilon(q_\varepsilon)),$$ where $\delta_2>0$ does not depend on $\varepsilon$. Let $\varepsilon\rightarrow0$ to get a contradiction.
\qed
\vspace{3mm}

In order to prove Theorem \ref{t3}, we need the solution of the Dirichlet problem for quaternionic Monge-Amp\`{e}re equation in quaternionic strictly pesudoconvex bounded domain in $\mathbb{H}^n$, which was solved by Alesker \cite{alesker4}.

Recall that an open bounded domain $\Omega\subseteq\mathbb{H}^n$ with a $C^\infty$-smooth boundary $\partial\Omega$ is called \emph{strictly pseudoconvex} if for every point $q_0\in \partial\Omega$ there exists a neighborhood $\mathcal {O}$ and a $C^\infty$-smooth strictly plurisubharmonic function $h$ on $\mathcal {O}$ such that $\Omega\cap\mathcal {O}=\{h<0\}$ and $\nabla h(q_0)\neq 0$.

\begin{lem}\label{l2.11}(Theorem 1.3 in \cite{alesker4}) Let $\Omega\subseteq\mathbb{H}^n$ be a bounded quaternionic strictly pseudoconvex domain. Let $0\leq f\in C(\overline{\Omega})$ and $\phi\in C(\partial \Omega)$. Then there exists a unique function $u$ such that
\equation\label{d1} \left\{
                      \begin{array}{ll}
                      u\in PSH(\Omega)\cap C(\overline{\Omega})&\\
                         det(u)=f ~~\text{in}~\Omega& \\
                        u|_{\partial\Omega}=\phi. &
                      \end{array}
                    \right.\endequation\end{lem}

\vspace{2mm}

 \emph{Proof of Theorem \ref{t3}} Denote by $\mathcal {F}$ the family of all subsolutions $\omega$ of (\ref{1}) satisfying $u\leq \omega\leq v$. By the viscosity comparison principle, $\mathcal {F}$ is non-empty. Let $$h:=\sup\{\omega:\omega\in \mathcal {F}\}.$$ By Choquet lemma, $h^*=(\limsup_j\omega_j)^*$, where $\omega_j$ is a sequence in $\mathcal {F}$. By Lemma \ref{l2.1}, $h^*$ is a subsolution of (\ref{1}).

We claim that $h_*$ is a supersolution of (\ref{1}). Assume by contradiction that $h_*$ is not a supersolution. Then by definition there exist $q_0\in\Omega$ and $\varphi\in C^2(\{q_0\})$ such that $\varphi$ touches $h_*$ from below at $q_0$, but
\begin{equation}\label{2.6}det( \varphi)(q_0)>F(q_0,\varphi(q_0)) .\end{equation}Then proceeding as in p.24 \cite{guide} we can construct a subsolution $\phi$ such that $\phi(q_j)>h(q_j)$ for some $q_j\in \Omega$. This contradiction leads to the conclusion that $h_*$ is a supersolution. For the reader's convenience, we briefly summarize the construction of $\phi$.

By the continuity of $F$ and (\ref{2.6}), we can find $r>0$ small enough such that $\varphi\leq h_*$ in a ball $B=B(q_0,r)$ and $det( \varphi)(q)>F(q,\varphi(q)),$ for all $q\in B$. Take $\varepsilon>0$ small enough and $0<\delta\ll\varepsilon$ such that $Q:=\varphi+\delta-\varepsilon|q-q_0|^2$ satisfies: $det( Q)(q)>F(q,Q(q)),$ for all $q\in B$. Define
$$\phi=\left\{
         \begin{array}{ll}
           \max\{h,Q\}, &\text{in}~B, \\
           h, & \text{in}~\Omega\backslash B.
         \end{array}
       \right.$$ Since $Q<h$ near $\partial B$, $\phi$ is plurisubharmonic and is a subsolution of (\ref{1}). Choose a sequence $\{q_j\}$ in $B$ converging to $q_0$ such that $h(q_j)\rightarrow h_*(q_0)$. Then $$Q(q_j)-h(q_j)\rightarrow Q(q_0)-h_*(q_0)=\varphi(q_0)+\delta-h_*(q_0)=\delta>0.$$
       It follows $\phi(q_j)=Q(q_j)>h(q_j)$, which contradicts the maximality of $h$.

Now we know that $h^*$ is a subsolution and $h_*$ is a supersolution. Since $g=u_*\leq h_*\leq h^*\leq v^*=g$ on $\partial\Omega$, by the viscosity comparison principle we get that $h=h_*=h^*$ is a continuous viscosity solution of (\ref{1}) with boundary value $g$. Now it suffices to show that $h$ is also a pluripotential solution of (\ref{1}).

It follows from Theorem \ref{t1} that $det( h)\geq F(q,h)$ in the pluripotential sense, i.e. $det(h)\geq F(q,h)$. Take an arbitrary ball $B$ in $\Omega$. Noting that $h$ and $F$ are continuous, by Lemma \ref{l2.11} we can solve the Dirichlet problem to find $\psi\in PSH(B)\cap C(\overline{B})$ with boundary value $h$ such that $det(\psi)=F(q,h)$ in $B$. By the pluripotential comparison principle (Proposition \ref{l2.10}) we have $h\leq \psi$ in $\overline{B}$. Let
$$\widetilde{\psi}=\left\{
                  \begin{array}{ll}
                    \psi, & \hbox{in}~B, \\
                    h, & \hbox{in}~\Omega\backslash\overline{B},
                  \end{array}
                \right. \text{and}\quad G(q)=F(q,h(q)).$$ Then $\widetilde{\psi}$ is a viscosity solution of $-det(\widetilde{\psi})+G(q)=0$. Again by viscosity comparison principle we have $\widetilde{\psi}\leq h$ in $\Omega$. Thus $\psi=h$ in $B$. The proof is complete.
\qed

\begin{cor} Let $\Omega$ be a quaternionic strictly pseudoconvex bounded domain in $\mathbb{H}^n$ and $g\in C(\partial\Omega)$. Assume that $F(q,t)$ is a continuous function on $\overline{\Omega}\times\mathbb{R}\rightarrow\mathbb{R}^+$ and is non-decreasing in the second variable. Then there exists a unique pluripotential (equiv. viscosity) solution to the Dirichlet problem (\ref{1}) with boundary value $g$.
\end{cor}
\proof It follows from \cite{alesker4} that there exists a continuous PSH function $u$ with boundary value $g$. Let $\rho\in C^2(\overline{\Omega})$ be a defining function of $\Omega$. Then $u+A\rho$ is a subsolution with boundary value $g$, provided $A$ is a sufficiently large constant. The smooth boundary implies the existence of harmonic function for arbitrary given continuous boundary data; it is a continuous supersolution to (\ref{1}). Therefore we obtain the existence and uniqueness by Theorem \ref{t3}.
\endproof
 \vskip 5mm

\section*{Acknowledgements}
This work is supported by Natural Science Foundation of SZU (grant no. 201424) and National Nature Science Foundation in China (No. 11401390; No.
11171298).

  \bibliographystyle{plain}
 \bibliography{mybibfile}





\end{document}